# Multifractal scaling of products of birth–death processes

VO V. ANH[1], NIKOLAI N. LEONENKO[2] and NARN-RUEIH SHIEH[3]

[1]*School of Mathematical Sciences, Queensland University of Technology, GPO Box 2434, Brisbane QLD 4001, Australia.* E-mail: *v.anh@qut.edu.au*
[2]*Cardiff School of Mathematics, Cardiff University, Senghennydd Road, Cardiff CF24 4AG, UK.* E-mail: *LeonenkoN@Cardiff.ac.uk*
[3]*Department of Mathematics, National Taiwan University, Taipei 10617, Taiwan.* E-mail: *shiehnr@math.ntu.edu.tw*

We investigate the scaling properties of products of the exponential of birth–death processes with certain given marginal discrete distributions and covariance structures. The conditions on the mean, variance and covariance functions of the resulting cumulative processes are interpreted in terms of the moment generating functions. We provide four illustrative examples of Poisson, Pascal, binomial and hypergeometric distributions. We establish the corresponding log-Poisson, log-Pascal, log-binomial and log-hypergeometric scenarios for the limiting processes, including their Rényi functions and dependence properties.

*Keywords:* geometric birth–death processes; log-binomial scenario; log-Pascal scenario; log-Poisson scenario; multifractal products

## 1. Introduction

Fully developed turbulence has been characterized by certain universal properties such as the scaling behavior

$$\mathrm{E}(v(x+l) - v(x))^q \sim l^{\zeta(q)}$$

across a distance $l$ of the $q$th-order moment of velocity fluctuations, or

$$\mathrm{E}(\varepsilon_l)^q \sim l^{\tau(q)}$$

for the $q$th-order moment of locally averaged energy dissipation $\varepsilon_l$ over a ball of size $l$. Kolmogorov's refined similarity hypothesis [20] leads to the relationship (Frisch [11])

$$\zeta(q) = \frac{q}{3} + \tau(q/3).$$

She and Lévêque [30] proposed a scaling model which predicts

$$\zeta(q) = \frac{q}{9} + 2\left[1 - \left(\frac{2}{3}\right)^{q/3}\right].$$

Dubrulle [8] explored further properties of the model and showed that its probability density function is related to the log-Poisson statistics of local non-dimensional energy dissipation. She and Waymire [31] proposed that the statistics of the inertial range of fully developed turbulence can be described by random multiplicative cascades, showed that these random cascades are log-Poisson and re-established the scaling model of She and Lévêque [30].

This paper will provide a construction of the log-Poisson random cascade and certain other models of the same type, such as log-Pascal, log-binomial and log-hypergeometric cascades, in the framework of multifractal products of birth–death processes. Multifractal products of stochastic processes were defined in Kahane [16, 17] and further investigated in Mannersalo *et al.* [22]. In the present paper, products of independent rescaled copies of a "mother" process are considered where the mother process is in the form of the exponential of a birth–death process with certain given marginal discrete distribution and covariance structure. The birth–death process is an important example of a homogeneous Markov chain $\{X(t),\ t \geq 0\}$ on the state space $S = \{-1, 0, 1, \ldots\}$. Here $X(t)$ can be considered as the size of a population at time $t$, which fluctuates according to the following rule: If at time $t$ the chain is in a state $i \in S$, in one transition it can go only to $i-1$ or $i+1$. The transition from state $i$ to state $i+1$ indicates a "birth" in the population, whereas the transition from $i$ to $i-1$ indicates a "death". A classification of birth–death processes was given in Karlin and McGregor [18]. A general theory of birth–death processes can be found in Dynkin and Yuskevich [9], Karlin and Taylor [19] and Bhattacharya and Waymire [6]. For the cases of Poisson, Pascal, binomial and hypergeometric distributions of these processes, we construct log-Poisson, log-Pascal, log-binomial and log-hypergeometric scenarios with non-trivial singularity spectra. Our method is based on an application of proper estimates of the maximal increments of the process in conjunction with the orthogonal expansion of its transition distribution.

There are many constructions of random multiplicative cascades ranging from the simple binomial cascade to measures generated by branching processes and the compound Poisson process (see Kahane [16, 17], Molchan [23], Falconer [10], Barral and Mandelbrot [5], Riedi [28], Mörters and Shieh [24, 25, 26] and Shieh and Taylor [32]). Many of these multifractal models were not designed to cover dependence structure and a natural form of the singularity spectrum (see Novikov [27] and Riedi [28], e.g.). Dependent cascades were first considered in Waymire and Williams [34], [35], including a general theory for Markov dependent cascade generators. Jaffard [15] showed that Lévy processes (except Brownian motion and Poisson process) are multifractal; but since the increments of a Lévy process are independent, this class excludes the effects of dependence. Moreover, Lévy processes have a linear singularity spectrum while real data often exhibit a strictly concave spectrum. In this paper, we pay attention to the strong correlation and nonlinear form of the singularity spectrum of a class of random cascades.

Our exposition relies on some results of Mannersalo *et al.* [22] on the basic properties of multifractal products of stochastic processes, but we provide a new interpretation of the conditions on the mean, variance and covariance functions of the resulting cumulative processes in terms of the moment generating functions. This interpretation is more useful for our development. We describe the behavior of the $q$th-order moments and Rényi functions, which are nonlinear, hence displaying the multifractality of the limiting cumulative processes. A property on their dependence structure, leading to their possible long-range dependence, is also obtained. The development of Mannersalo *et al.* [22] can be applied to a large class of stationary processes such as stationary diffusion processes. As applications, Mannersalo *et al.* [22] looked at continuous-time two-state Markov processes, while in the present paper we also consider infinite state space and some new scenarios such as those cited above.

## 2. Multifractal products of stochastic processes

This section recaptures some basic results on multifractal products of stochastic processes as developed in Kahane [16, 17] and Mannersalo *et al.* [22]. We provide a new interpretation of their conditions based on the moment generating functions, which is useful for our exposition.

We introduce the following conditions:

**A′**. Let $\Lambda^{(i)}(t)$, $t \in [0,1]$, $i = 1, 2, \ldots$, be a sequence of independent, strictly stationary, positive stochastic processes such that, for all $t, t_1, t_2 \in [0,1]$ and $i = 0, 1, 2, \ldots$, the following assumptions hold:

$$E\Lambda^{(i)}(t) = 1, \tag{2.1}$$

$$\operatorname{Var}\Lambda^{(i)}(t) = \sigma_\Lambda^2 < \infty, \tag{2.2}$$

$$\operatorname{Cov}(\Lambda^{(i)}(t_1), \Lambda^{(i)}(t_2)) = R_\Lambda(t_1 - t_2) = \sigma_\Lambda^2 \rho_i(t_1 - t_2), \qquad \rho_i(0) = 1. \tag{2.3}$$

We consider the following setting:

**A″**. Let $\Lambda_b^{(i)}$ be independent rescaled copies of some measurable separable mother process $\Lambda$, that is,

$$\Lambda_b^{(i)}(t) \stackrel{d}{=} \Lambda^{(i)}(tb^i), \qquad t \in [0,1],\ i = 0, 1, 2, \ldots,$$

where the scaling parameter $b > 1$, $E\Lambda(t) = 1$, and $\stackrel{d}{=}$ denotes equality in finite-dimensional distributions.

Moreover, in the examples of Section 4, the stationary mother process satisfies the following conditions:

**A‴**. For $t \in [0,1]$, let $\Lambda(t) = \exp\{X(t)\}$, where $X(t)$ is a stationary process with $EX^2(t) < \infty$,

$$\operatorname{Cov}(X(t_1), X(t_2)) = R_X(t_1 - t_2) = \sigma_X^2 r_X(t_1 - t_2), \qquad r_X(0) = 1.$$

We assume that there exist a marginal probability mass function $p_\theta(x)$ and a bivariate probability mass function $p_\theta(x_1, x_2; t_1 - t_2)$ such that the moment generating function

$$M(\zeta) = \mathrm{E}\exp\{\zeta X(t)\}$$

and the bivariate moment generating function

$$M(\zeta_1, \zeta_2; t_1 - t_2) = \mathrm{E}\exp\{\zeta_1 X(t_1) + \zeta_2 X(t_2)\}$$

exist. Here, $\theta$ is the parameter vector of the mass function of the process $X(t)$.

Under the conditions $\mathbf{A}'$–$\mathbf{A}'''$, the assumptions (2.1)–(2.3) can be rewritten as

$$\mathrm{E}\Lambda_b^{(i)}(t) = M(1) = 1;$$

$$\mathrm{Var}\,\Lambda_b^{(i)}(t) = M(2) - 1 = \sigma_\Lambda^2 < \infty;$$

$$\mathrm{Cov}(\Lambda_b^{(i)}(t_1), \Lambda_b^{(i)}(t_2)) = M(1, 1; (t_1 - t_2)b^i) - 1, \ b > 1.$$

We define the finite product processes

$$\Lambda_n(t) = \prod_{i=0}^n \Lambda_b^{(i)}(t) = \exp\left\{\sum_{i=0}^n X^{(i)}(tb^i)\right\},$$

and the cumulative processes

$$A_n(t) = \int_0^t \Lambda_n(s)\,\mathrm{d}s, \qquad n = 0, 1, 2, \ldots,$$

where $X^{(i)}(t), i = 0, \ldots, n, \ldots$, are independent copies of a stationary process $X(t), t \geq 0$.

We also consider the corresponding positive random measures defined on Borel sets $B$ of $[0,1]$:

$$\mu_n(B) = \int_B \Lambda_n(s)\,\mathrm{d}s, \qquad n = 0, 1, 2, \ldots.$$

Kahane [17] proved that the sequence of random measures $\mu_n$ converges weakly almost surely to a random measure $\mu$. Moreover, given a finite or countable family of Borel sets $B_j$ on $[0,1]$, it holds that $\lim_{n\to\infty} \mu_n(B_j) = \mu(B_j)$ for all $j$ with probability one. The almost sure convergence of $A_n(t)$ in countably many points of $[0, 1]$ can be extended to all points in $[0, 1]$ if the limit process $A(t)$ is almost surely continuous. In this case, $\lim_{n\to\infty} A_n(t) = A(t)$ with probability one for all $t \in [0, 1]$. As noted in Kahane [17], there are two extreme cases: (i) $A_n(t) \to A(t)$ in $L_1$ for each given $t$, in which case $A(t)$ is not almost surely zero and is said to be fully active (non-degenerate) on $[0, 1]$; (ii) $A_n(1)$ converges to 0 almost surely, in which case $A(t)$ is said to be degenerate on $[0, 1]$. Sufficient conditions for non-degeneracy and degeneracy in a general situation and relevant examples are provided in Kahane [17] (equations (18) and (19), respectively.) The condition for complete degeneracy is detailed in Theorem 3 of [17].

The Rényi function, also known as the deterministic partition function, is defined as

$$T(q) = \liminf_{n\to\infty} \frac{\log \mathrm{E} \sum_{k=0}^{2^n-1} \mu^q(I_k^{(n)})}{\log |I_k^{(n)}|}$$

$$= \liminf_{n\to\infty} \left(-\frac{1}{n}\right) \log_2 \mathrm{E} \sum_{k=0}^{2^n-1} \mu^q(I_k^{(n)}),$$

where $I_k^{(n)} = [k2^{-n}, (k+1)2^{-n}]$, $k = 0, 1, \ldots, 2^n - 1$, $|I_k^{(n)}|$ is its length and $\log_b$ is log to the base $b$.

**Remark 1.** The multifractal formalism for random cascades can be stated in terms of the Legendre transform of the Rényi function:

$$T^*(\alpha) = \min_{q\in\mathbb{R}}(q\alpha - T(q)).$$

In fact, let $f(\alpha)$ be the Hausdorff dimension of the set

$$C_\alpha = \left\{ t \in [0,1] : \lim_{n\to\infty} \frac{\log \mu(I_k^{(n)}(t))}{\log |I_k^{(n)}|} = \alpha \right\},$$

where $I_k^{(n)}(t)$ is a sequence of intervals $I_k^{(n)}$ that contain $t$. The function $f(\alpha)$ is known as the singularity spectrum of the measure $\mu$, and we refer to $\mu$ as a multifractal measure if $f(\alpha) \neq 0$ for a continuum of $\alpha$ (Lau [21]). In order to determine the function $f(\alpha)$, Hentschel and Procaccia [14], Frisch and Parisi [12] and Halsey *et al.* [13], for example, proposed to use the relationship

$$f(\alpha) = T^*(\alpha). \tag{2.4}$$

This relationship may not hold for a given measure (see, e.g., Taylor [33]). When the equality (2.4) is established for a measure $\mu$, we say that the multifractal formalism holds for this measure.

Mannersalo *et al.* [22] presented the conditions for $L_2$-convergence and scaling of moments.

**Theorem 1** *(Mannersalo, Norros and Riedi [22])*. *Suppose that the conditions* $\mathbf{A}'$–$\mathbf{A}'''$ *hold.*

*If, for some positive numbers $\delta$ and $\gamma$,*

$$\exp\{-\delta|\tau|\} \leq \rho(\tau) = \frac{M(1,1;\tau) - 1}{M(2) - 1} \leq |C\tau|^{-\gamma}, \tag{2.5}$$

*then $A_n(t)$ converges in $L_2$ if and only if*

$$b > 1 + \sigma_\Lambda^2 = M(2).$$

If $A_n(t)$ converges in $L_1$, then the limit process $A(t)$ satisfies the recursion

$$A(t) = \frac{1}{b} \int_0^t \Lambda(s) \, d\tilde{A}(bs), \qquad (2.6)$$

where the processes $\Lambda(t)$ and $\tilde{A}(t)$ are independent, and the processes $A(t)$ and $\tilde{A}(t)$ have identical finite-dimensional distributions.

If $A(t)$ is non-degenerate, the recursion (2.6) holds, $A(1) \in L_q$ for some $q > 0$ and $\sum_{n=0}^{\infty} c(q, b^{-n}) < \infty$, where $c(q,t) = \operatorname{E}\sup_{s \in [0,t]} |\Lambda^q(0) - \Lambda^q(s)|$, then there exist constants $\overline{C}$ and $\underline{C}$ such that, for all $t \in [0,1]$,

$$\underline{C} t^{q - \log_b \operatorname{E}\Lambda^q(t)} \leq \operatorname{E}A^q(t) \leq \overline{C} t^{q - \log_b \operatorname{E}\Lambda^q(t)}, \qquad (2.7)$$

which will be written as

$$\operatorname{E}A^q(t) \sim t^{q - \log_b \operatorname{E}\Lambda^q(t)}.$$

If, on the other hand, $A(1) \in L_q$, $q > 1$, then the Rényi function is given by

$$T(q) = q - 1 - \log_b \operatorname{E}\Lambda^q(t) = q - 1 - \log_b M(q). \qquad (2.8)$$

If $A(t)$ is non-degenerate, $A(1) \in L_2$, and $\Lambda(t)$ is positively correlated, then

$$\operatorname{Var} A(t) \geq \operatorname{Var} \int_0^t \Lambda(s) \, ds. \qquad (2.9)$$

Hence, if $\int_0^t \Lambda(s) \, ds$ is strongly dependent, then $A(t)$ is also strongly dependent.

**Remark 2.** The result (2.7) means that the process $A(t)$, $t \in [0,1]$ with stationary increments behaves as

$$\log \operatorname{E}[A(t+\delta) - A(t)]^q \sim K(q) \log \delta + C_q \qquad (2.10)$$

for a wide range of resolutions $\delta$ with a nonlinear function

$$K(q) = q - \log_b \operatorname{E}\Lambda^q(t) = q - \log_b M(q),$$

where $C_q$ is a constant. In this sense, the stochastic process $A(t)$ is said to be multifractal. The function $K(q)$, which contains the scaling parameter $b$ and all the parameters of the marginal distribution of the mother process $X(t)$, can be estimated by running the regression (2.10) for a range of values of $q$. For the four examples in Section 4, the explicit form of $K(q)$ is obtained. Hence these parameters can be estimated by minimizing the mean square error between the $K(q)$ curve estimated from data and its analytical form for a range of values of $q$. This method has been used for multifractal characterization of complete genomes in Anh *et al.* [2].

## 3. Geometric birth–death processes

We consider the one-dimensional birth–death process $X = \{X(t), \ t \geq 0\}$ on the state space $S = \{-1, 0, 1, \ldots\}$. This process is defined by its birth–death rates $\lambda_i, \mu_i, i \in S$ (Karlin and McGregor [18]), that is, $X$ is a Markov process with stationary transition probabilities

$$P_{ij}(t) = P\{X(t+s) = j | X(s) = i\}, \qquad i, j \in S,$$

which are independent of $s$. In addition, we assume that $P_{ij}(t)$ satisfy

$$P_{i,i+1}(h) = \lambda_i h + \mathrm{o}(h), \qquad h \to 0, \ i \in S; \tag{3.1}$$

$$P_{i,i-1}(h) = \mu_i h + \mathrm{o}(h), \qquad h \to 0, \ i \in S; \tag{3.2}$$

$$P_{i,i}(h) = 1 - (\lambda_i + \mu_i)h + \mathrm{o}(h), \qquad h \to 0, \ i \in S; \tag{3.3}$$

$$P_{i,j}(0) = \delta_i^j; \qquad P_{-1,-1}(t) = 1, \qquad P_{-1,i}(t) = 0, \qquad t \geq 0, \ i \neq -1.$$

The order $\mathrm{o}(h)$ may depend on $i$, and thus we will use the notation $\mathrm{o}(h; i)$ when it becomes necessary. We assume $\mu_0 \geq 0, \lambda_0 > 0, \lambda_i, \mu_i > 0, i \geq 1$.

If $\mu_0 > 0$, then we have an absorbing state $-1$; once the process enters $-1$, it can never leave it. If $\mu_0 = 0$, we have a reflecting state $0$. After entering $0$, the process will always go back to state $1$ after some time. In this case, the state $-1$ can never be reached and so can be ignored, and we take $S = \{0, 1, 2, \ldots\}$.

We may then suppose that $\mu_0 = 0$, so that the state $-1$ is ignored. We assume that

$$\pi := \sum_{k=0}^{\infty} \pi_k < \infty, \qquad \sum_{k=0}^{\infty} \frac{1}{\lambda_k \pi_k} = \infty,$$

where the potential coefficients $\pi_i$ are given by

$$\pi_i = \frac{\lambda_0 \lambda_1 \cdots \lambda_{i-1}}{\mu_1 \mu_2 \cdots \mu_i}, \qquad i \geq 1, \pi_0 = 1,$$

so that the process is ergodic (see Karlin and Magregor [18], Theorem 2). Then,

$$\lim_{t \to \infty} P_{i,j}(t) := p_j = \frac{\pi_j}{\pi}, \qquad j \geq 0,$$

exist and are independent of the initial state $i$. We write the stationary distribution as $P = (p_j, \ j = 0, 1, 2, \ldots)$. The infinitesimal generator $\mathcal{A}$ of the process is given by

$$\mathcal{A}f(i) = \lambda_i f(i+1) - (\lambda_i + \mu_i)f(i) + \mu_i f(i-1)$$

for all bounded real-valued functions $f$.

We denote by $l_2(\mathbb{N}, P)$ the Hilbert space of functions $f(j), j = 0, 1, \ldots$, such that

$$\sum_{j=0}^{\infty} f^2(j) p_j < \infty.$$

The spectral analysis of birth–death processes is based on the birth–death polynomials $Q_n(x), n = 0, 1, 2, \ldots$, which are defined by the recursive relations

$$-xQ_n(x) = \mu_n Q_{n-1}(x) - (\lambda_n + \mu_n)Q_n(x) + \lambda_n Q_{n+1}(x), \qquad n = 0, 1, 2, \ldots, \qquad (3.4)$$

with $Q_{-1}(x) = 0, Q_0(x) = 1$. Karlin and McGregor [18] proved that there exists a positive Borel measure $\phi(\mathrm{d}x)$ with total mass 1 and with support on the non-negative real numbers, called the spectral measure of the process, such that the transition probabilities are represented in the spectral form as

$$P_{i,j}(t) = \pi_j \int_0^\infty \mathrm{e}^{-xt} Q_i(x) Q_j(x) \phi(\mathrm{d}x).$$

The polynomials $Q_n(x)$ are orthogonal with respect to $\phi$. For certain choices of birth–death parameters, we can apply the theory of classical orthogonal polynomials with respect to discrete distributions. In the case $\mu_0 = 0$, which we assume in this paper, we have $Q_n(0) = 1, n \geq 0$ and that 0 is $\phi$-atomic; the stationary distribution $P$ is then given by $p_j = \pi_j \phi(\{0\})$.

The following proposition plays a key role in our exposition. We note that the assumptions (ii) and (iii) in Section 4 of Mannersalo *et al.* [22] generally do not hold in the case that the state space $S$ is infinite; indeed their assumption (ii) assumes that the jump rate $\nu(x)$ is uniformly bounded from below and from above for all $x$; this assumption is not suitable for an infinite state space, as our present Section 4 shows. We also remark that the condition (3.5) below, which is needed for the validity of Proposition 1 for an infinite state space, may not be the best, and the possibility of improving it is left for a future study (we thank the referee for noticing this). In the following, $\mathrm{E}^P$ denotes expectation with respect to the stationary distribution $P$.

**Proposition 1.** *We assume, besides the strict stationarity, that the small order of magnitude in the transition probabilities (3.1)–(3.3) is such that $\mathrm{o}(h; k) = k^\delta \mathrm{o}(h; 1)$ for all $k = 1, 2, 3, \ldots$, for some $\delta > 0$. Moreover, we assume that*

$$S_q := \sum_{k=0}^\infty \pi_k (\mathrm{e}^{qk}(\lambda_k + \mu_k) + k^\delta) < \infty, \qquad q \in \mathbb{R}. \qquad (3.5)$$

*Let $g(x) := \mathrm{e}^{q(x-c)}$, where $c \in \mathbb{R}$. Then, for any $b > 1$,*

$$\sum_{n=0}^\infty \mathrm{E}^P \left( \max_{s \in [0, b^{-n}]} |g(X(s)) - g(X(0))| \right) < \infty. \qquad (3.6)$$

**Proof.** We first note that

$$\mathrm{E}^P[\cdot] = \frac{1}{\pi} \sum_k \pi_k \mathrm{E}(\cdot | X(0) = k). \qquad (3.7)$$

Since the process is cadlag and integer-valued, we have, for each $t > 0$,

$$E\left(\max_{0<s\leq t} |g(X(s)) - g(X(0))| \mid X(0) = k\right)$$
$$\leq E\left(\sum_{0<s\leq t, \Delta X_s \neq 0} |g(X(s)) - g(X(s-))| \;\Big|\; X(0) = k\right),$$

where $\Delta X_s = X(s) - X(s-), 0 < s \leq t$, and those $s$ for which $\Delta X_s \neq 0$ are countably many. The defining property of birth–death processes implies that

$$P\{\Delta X_s \neq 0 | X(s-)\} \leq const \times (\lambda_{X(s-)} + \mu_{X(s-)}),$$

where the constant is taken to be the expected jump time $+1$, which is independent of $s$. We note that, by the definition of $g(\cdot)$,

$$0 < |g(X(s)) - g(X(s-))| \leq e^{qh} \max\{|e^q - 1|, |1 - e^{-q}|\}$$

when $X(s-) = h$ and $\Delta X_s \neq 0$. As $t \downarrow 0$, the joint probabilities of $\{X(s-) = h, X(0) = k\}$, $h$ varying over $S - \{k\}$, are all of order $o(t; k)$, except when $h = k+1$ and $h = k-1$, in which case they are of order $\lambda_k t$ and $\mu_k t$, respectively. Therefore, for large $n$ and by the strong Markov property of the process, which starts afresh at each jump time, we have

$$E\left(\max_{0<s\leq b^{-n}} |g(X(s)) - g(X(0))| \mid X(0) = k\right)$$
$$\leq E\left(\sum_{0<s\leq t, \Delta X_s \neq 0} |g(X(s)) - g(X(s-))| \;\Big|\; X(0) = k\right)$$
$$\leq E\left(\sum_{0<s\leq t, \Delta X_s \neq 0} E(|g(X(s)) - g(X(s-))| \mid X(s-)) \;\Big|\; X(0) = k\right)$$
$$\leq E_{X(0)=k}\left(\sum_{0<s\leq b^{-n}, \Delta X_s \neq 0} E(|g(X(s)) - g(X(s-))| \;\Big|\; X(s-) = k \pm 1)\right)$$
$$+ o(b^{-n}; k)$$
$$\leq C_q(e^{qk}(\lambda_k + \mu_k) + k^\delta)b^{-n}.$$

In the above inequalities, we have arranged, without indexing explicitly, the jump times in successively increasing order, and have made use of the assumption $o(h; k) = k^\delta o(h; 1)$ for all $k = 1, 2, 3, \ldots$. We note that the expected number of the jumps is at most 1 on $[0, b^{-n}]$, and thus the summation of jumps in the above is indeed reduced to one term only. In view of the conditional expectation formula (3.7) and the definition of $S_q$, which we assume to be finite, we get

$$\sum_{n=0}^{\infty} \mathrm{E}^P \left( \max_{s \in [0, b^{-n}]} |g(X(s)) - g(X(0))| \right) \leq (C_q \cdot S_q) \sum_{n=0}^{\infty} b^{-n}.$$

This completes the proof of the proposition. □

**Remark 3.** We usually consider $q > 0$, in which case the finiteness of $S_q$ depends on the growth of the sequences $\lambda_k, \mu_k$. We have included here the case $q \leq 0$ since it incurs no additional difficulty.

## 4. Multifractal scenarios

This section introduces four illustrative multifractal scenarios. The mother process will take the form $\Lambda(t) = \exp\{X(t) - c_X\}$, where $X(t)$ is a stationary birth–death process and $c_X$ is a constant depending on the parameters of its marginal distribution. This form is needed for the condition $\mathrm{E}\Lambda(t) = 1$ to hold.

### 4.1. Log-Poisson scenario

The log-Poisson statistics in fully developed turbulence were discussed in Dubrulle [8] and She and Waymire [31]. In this section, we provide a related model, namely the log-Poisson scenario for multifractal products of stochastic processes.

**B′.** Consider a mother process of the form

$$\Lambda(t) = \exp\left\{X(t) - \frac{\lambda}{\mu}(\mathrm{e} - 1)\right\}, \tag{4.1}$$

where $X(t)$, $t \geq 0$ is a stationary birth–death process with marginal Poisson distribution $\mathrm{Poi}(\frac{\lambda}{\mu}), \lambda > 0, \mu > 0$, and rates

$$\lambda_n = \lambda, \qquad \mu_n = \mu n, \qquad n \geq 0.$$

The covariance function of the process $X(t)$ then takes the form

$$R_X(t) = \frac{\lambda}{\mu} r_X(t), \qquad r_X(t) = \mathrm{e}^{-\mu t}.$$

Under condition **B′**, we obtain the following moment generating function:

$$M(\zeta) = \mathrm{E} \exp\left\{\zeta\left(X(t) - \frac{\lambda}{\mu}(\mathrm{e} - 1)\right)\right\}$$

$$= \exp\left\{-\zeta\frac{\lambda}{\mu}(\mathrm{e} - 1) + \frac{\lambda}{\mu}(\mathrm{e}^\zeta - 1)\right\}, \qquad \zeta \in \mathbb{R}.$$

It turns out that, in this case

$$\log_b \mathrm{E}\Lambda^q(t) = \frac{-q\lambda/\mu(\mathrm{e}-1) + \lambda/\mu(\mathrm{e}^q-1)}{\log b}, \qquad q > 0. \tag{4.2}$$

We can formulate the following:

**Theorem 2.** *Suppose that condition* **B'** *holds. Then, for any*

$$b > \exp\left\{\frac{\lambda}{\mu}(\mathrm{e}-1)^2\right\},$$

*the stochastic processes*

$$A_n(t) = \int_0^t \prod_{j=0}^n \Lambda^{(j)}(sb^j)\,\mathrm{d}s, \qquad t \in [0,1],$$

*converge in $L_2$ to the stochastic process $A(t), t \in [0,1]$ as $n \to \infty$ such that, if $A(1) \in L_q$ for $q \in (0, \infty)$,*

$$\mathrm{E}A^q(t) \sim t^{q(1+\lambda/(\mu \log b)(\mathrm{e}-1)) - (1/\log b)(\lambda/\mu)(\mathrm{e}^q-1)},$$

*and the Rényi function is given by*

$$T(q) = q\left(1 + \frac{\lambda}{\mu \log b}(\mathrm{e}-1)\right) - \frac{1}{\log b}\frac{\lambda}{\mu}(\mathrm{e}^q-1) - 1.$$

*Moreover,*

$$\mathrm{Var}\,A(t) \geq \int_0^t \int_0^t (\mathrm{e}^{\lambda/\mu \mathrm{e}^{-\mu|u-v|}} - 1)\,\mathrm{d}u\,\mathrm{d}v \geq 2t\frac{\lambda}{\mu}\frac{(1-\mathrm{e}^{-\mu t}) + 1 - \mathrm{e}^{-\mu t}(\mu t+1)}{\mu^2}.$$

**Proof.** We have $\delta = 1$, and

$$\pi_k \leq \mathit{const} \times \frac{1}{k!}\left(\frac{\lambda}{\mu}\right)^k.$$

Thus,

$$S_q := \sum_{k=0}^\infty \pi_k (\mathrm{e}^{qk}(\lambda_k + \mu_k) + k)$$

$$\leq \mathit{const} \times \left(\lambda \sum_{k=0}^\infty \frac{1}{k!}\left(\frac{\lambda}{\mu}\mathrm{e}^q\right)^k + \mu \sum_{k=1}^\infty \frac{1}{(k-1)!}\left(\frac{\lambda}{\mu}\mathrm{e}^q\right)^k + \sum_{k=1}^\infty \frac{1}{(k-1)!}\left(\frac{\lambda}{\mu}\right)^k\right) < \infty$$

for all positive $\lambda, \mu, q$. The condition for the inequalities (2.7) of Theorem 1 to hold then follows from Proposition 1. Now we consider the correlation decay and show that the condition (2.5) of Theorem 1 holds. To see this, we consider the mother process

$$\Lambda(t) = G(X(t)), \qquad G(u) = \exp\left\{u - \frac{\lambda}{\mu}(\mathrm{e} - 1)\right\}$$

as a nonlinear transformation of the Markov process $X(t), t \in [0,1]$ which has marginal distribution

$$p_j = P(X(t) = j) = \frac{\mathrm{e}^{-\lambda/\mu}(\lambda/\mu)^j}{j!}, \qquad j = 0, 1, 2, \ldots.$$

Let $C_n(j; \frac{\lambda}{\mu})$, $n = 0, 1, 2, \ldots$, be Charlier polynomials, defined via the generating function

$$\sum_{n=0}^{\infty} C_n\left(j; \frac{\lambda}{\mu}\right) \frac{z^n}{n!} = \mathrm{e}^z \left(1 - \frac{z\mu}{\lambda}\right)^j, \qquad 0 < z < \frac{\mu}{\lambda}$$

(see, e.g., Chihara [7] and Schoutens [29]). These polynomials form a complete system of orthogonal polynomials in the Hilbert space $l_2(\mathbb{N}, \mathrm{Poi}(\frac{\lambda}{\mu}))$. In this case, the recursive equation (3.4) for birth–death polynomials $Q_n(x) = C_n(\frac{x}{\mu}, \frac{\lambda}{\mu})$ is the recursive equation for the Charlier polynomials:

$$0 = nC_{n-1}\left(x, \frac{\lambda}{\mu}\right) + \left(x - \frac{\lambda}{\mu} - n\right)C_n\left(x, \frac{\lambda}{\mu}\right) + \frac{\lambda}{\mu}C_{n+1}\left(x, \frac{\lambda}{\mu}\right), \qquad n \geq 0,$$

where

$$C_0\left(x, \frac{\lambda}{\mu}\right) = 1, \qquad C_{-1}\left(x, \frac{\lambda}{\mu}\right) = 0.$$

Then, the following expansion of the bivariate distribution holds:

$$P(X(t) = j, X(s) = k) = p_j p_k \sum_{n=0}^{\infty} \mathrm{e}^{-\mu n |t-s|} C_n\left(j; \frac{\lambda}{\mu}\right) C_n\left(k; \frac{\lambda}{\mu}\right) \frac{(\lambda/\mu)^n}{n!}.$$

Note that

$$G(u) \in l_2\left(\mathbb{N}, \mathrm{Poi}\left(\frac{\lambda}{\mu}\right)\right),$$

since

$$\sum_{j=0}^{\infty} \mathrm{e}^{2(j - (\lambda/\mu)(\mathrm{e}-1))} \mathrm{e}^{-\lambda/\mu} \frac{(\lambda/\mu)^j}{j!} = \exp\left\{-\frac{\lambda}{\mu}(\mathrm{e}-1)^2\right\} < \infty.$$

We also note that $G(u)$ has the expansion

$$G(u) = \sum_{k=0}^{\infty} a_k C_k\left(u; \frac{\lambda}{\mu}\right) \frac{1}{d_k}, \qquad d_k^2 = \frac{k! \mu^k}{\lambda^k},$$

$$a_k = \sum_{n=0}^{\infty} C_k\left(n; \frac{\lambda}{\mu}\right) e^{-\lambda/\mu} \frac{(\lambda/\mu)^n}{n!} e^{n-(\lambda/\mu)(e-1)} \frac{1}{d_k},$$

$$\sum_{k=0}^{\infty} \frac{a_k^2}{d_k^2} = \sum_{k=0}^{\infty} e^{-\lambda/\mu} \frac{(\lambda/\mu)^k}{k!} e^{2(k-(\lambda/\mu)(e-1))} < \infty.$$

Taking into account the orthogonality

$$\sum_{j=0}^{\infty} C_n\left(j; \frac{\lambda}{\mu}\right) C_m\left(j; \frac{\lambda}{\mu}\right) e^{-\lambda/\mu} \frac{(\lambda/\mu)^j}{j!} = \delta_n^m d_n^2, \qquad n, m \geq 1,$$

where $\delta_k^j$ is the Kronecker delta function, we obtain

$$\Lambda(t) = \sum_{k=0}^{\infty} \frac{a_k}{d_k} C_k\left(X(t); \frac{\lambda}{\mu}\right).$$

Thus,

$$R_\Lambda(t,s) = \mathrm{Cov}(\Lambda(t), \Lambda(s)) = \sum_{k=1}^{\infty} \frac{a_k^2}{d_k^2} e^{-k\mu|t|},$$

and

$$e^{-\mu|\tau|} \frac{a_1^2}{d_1^2} \leq R_\Lambda(\tau) \leq e^{-\mu|\tau|} \sum_{k \geq 1} \frac{C_k^2}{d_k^2}.$$

This completes the proof. □

### 4.2. Log-Pascal scenario

**B″**. Consider a mother process of the form

$$\Lambda(t) = \exp\{X(t) - c_X\}, \qquad c_X = \beta \log \frac{1 - \lambda/\mu}{1 - (1 - \lambda/\nu)e^{\mu-\lambda}}, \tag{4.3}$$

where $X(t)$, $t \geq 0$ is a stationary birth–death process with marginal Pascal-type distribution $\mathrm{Pas}(\beta, \frac{\lambda}{\mu})$:

$$p_j = P(X(0) = (\mu-\lambda)j) = \left(1 - \frac{\lambda}{\mu}\right)^\beta \frac{(\beta)_j}{j!} \left(\frac{\lambda}{\mu}\right)^j, \qquad j = 0, 1, 2, \ldots, \ 0 < \lambda < \mu, \ \beta > 0,$$

and rates

$$\lambda_n = (n+\beta)\lambda, \qquad \mu_n = \mu n.$$

In the above expression for $p_j$, we have used the Pochhammer symbol $(\beta)_j := \beta \cdot (\beta + 1) \cdots (\beta + j - 1), j \geq 1, (\beta)_0 = 1$. We note that

$$\pi_j = \frac{(\beta)_j}{j!} \left(\frac{\lambda}{\mu}\right)^j, \qquad j = 0, 1, 2, \ldots.$$

Therefore the sum $S_q$ in Proposition 1 is finite when $q < \log(\frac{\mu}{\lambda})$.

Note that

$$EX(t) = (\mu - \lambda)\beta \frac{\lambda}{\mu} \bigg/ \left(1 - \frac{\lambda}{\mu}\right), \qquad \text{Var}\, X(t) = (\mu - \lambda)^2 \beta \frac{\lambda}{\mu} \bigg/ \left(1 - \frac{\lambda}{\mu}\right)^2,$$

$$R_X(t) = \text{Var}\, X(t) r_X(t), \qquad r_X(t) = e^{-\mu t}, \; t \geq 0.$$

Under condition **B″**, we obtain the following moment generating function:

$$M(\zeta) = E \exp\{\zeta X(t)\} = \left(\frac{1 - \lambda/\mu}{1 - (1 - \lambda/\mu)e^{\zeta(\mu-\lambda)}}\right)^\beta, \qquad 0 < \zeta < -\log\frac{\lambda}{\mu}\bigg/(\mu - \lambda).$$

We can formulate the following:

**Theorem 3.** *Suppose that condition* **B″** *holds. Then, for any*

$$b > e^{-2c_X} \left(\frac{1 - \lambda/\mu}{1 - (1 - \lambda/\mu)e^{2(\mu-\lambda)}}\right)^\beta,$$

*the stochastic processes*

$$A_n(t) = \int_0^t \prod_{j=0}^n \Lambda^{(j)}(sb^j)\, ds, \qquad t \in [0,1],$$

*converge in $L_2$ to the stochastic process $A(t)$, $t \in [0,1]$ as $n \to \infty$ such that, if $A(1) \in L_q$ for $q \in (0, \min\{\log \frac{\mu}{\lambda}, -\log \frac{\lambda}{\mu}/(\mu - \lambda)\})$,*

$$EA^q(t) \sim t^{K(q)},$$

*where*

$$K(q) = q\left(1 + \frac{c_X}{\log b}\right) + \frac{1}{\log b}\log\left(1 - \left(1 - \frac{\lambda}{\mu}\right)e^{q(\mu-\lambda)}\right) - \beta\frac{\log(1 - \lambda/\mu)}{\log b},$$

*and the Rényi function is given by*

$$T(q) = q\left(1 + \frac{c_X}{\log b}\right) + \frac{1}{\log b}\log\left(1 - \left(1 - \frac{\lambda}{\mu}\right)e^{q(\mu-\lambda)}\right) - \beta\frac{\log(1 - \lambda/\mu)}{\log b} - 1.$$

*Moreover,*

$$\operatorname{Var} A(t) \geq \int_0^t \int_0^t (\mathrm{e}^{\operatorname{Var} X(t) \mathrm{e}^{-\mu|u-v|}} - 1) \, \mathrm{d}u \, \mathrm{d}v$$

$$\geq 2t(\mu - \lambda)^2 \beta \frac{\lambda}{\mu} \frac{(1 - \mathrm{e}^{-\mu t}) + 1 - \mathrm{e}^{-\mu t}(\mu t + 1)}{\mu^2}.$$

**Proof.** We have again $\delta = 1$. In view of $q < \log(\frac{\mu}{\lambda})$ and

$$S_q \leq \text{const} \times \sum_k k \left( \frac{\lambda}{\mu} \mathrm{e}^q \right)^k,$$

the condition for the inequalities (2.7) of Theorem 1 to hold then follows from Proposition 1. We now show that the condition (2.5) of Theorem 1 holds. For this purpose, we consider the mother process

$$\Lambda(t) = G(X(t)), \qquad G(u) = \exp\{u - c_X\}$$

as a nonlinear transformation of the Markov process $X(t), t \in [0, 1]$, which has marginal distribution $p_j \sim \operatorname{Pas}(\beta, \frac{\lambda}{\mu})$ and the following bivariate distribution expansion:

$$P(X(t) = j, X(s) = k) = p_j p_k \sum_{n=0}^{\infty} \frac{\mathrm{e}^{-\mu n |t-s|}}{d_n} M_n\left(j; \beta, \frac{\lambda}{\mu}\right) M_n\left(k; \beta, \frac{\lambda}{\mu}\right),$$

$$d_n^2 = \frac{n!}{(\lambda/\mu)^n (\beta)_n},$$

where $M_n(j; \beta, \frac{\lambda}{\mu}), n = 0, 1, \ldots$ are Meixner polynomials, defined via the generating function

$$\sum_{n=0}^{\infty} (\beta)_n M_n\left(j; \beta, \frac{\lambda}{\mu}\right) \frac{z^n}{n!} = (1-z)^{-j-\beta}\left(1 - \frac{z\mu}{\lambda}\right)^j, \qquad 0 < z < \frac{\lambda}{\mu}$$

(see, e.g., Chihara [7] and Schoutens [29]). The recursive equation for the birth–death polynomial (3.4) becomes in this case the recursive equation for the Meixner polynomials:

$$-jM_n\left(j; \beta, \frac{\lambda}{\mu}\right) = b_n M_{n+1}\left(j; \beta, \frac{\lambda}{\mu}\right) + \gamma_n M_n\left(j; \beta, \frac{\lambda}{\mu}\right) + c_n M_{n-1}\left(j; \beta, \frac{\lambda}{\mu}\right), \qquad n \geq 0,$$

where

$$M_{-1}\left(j; \beta, \frac{\lambda}{\mu}\right) = 0, \qquad M_0\left(j; \beta, \frac{\lambda}{\mu}\right) = 1,$$

and

$$b_n = \frac{\lambda}{\mu}(n+\beta)\frac{1}{1 - \lambda/\mu}, \qquad \gamma_n = -\frac{n + (\lambda/\mu)(n+\beta)}{1 - \lambda/\mu}, \qquad c_n = \frac{n}{1 - \lambda/\mu}.$$

We note that
$$G(u) \in l_2\left(\mathbb{N}, \mathrm{Pas}\left(\beta, \frac{\lambda}{\mu}\right)\right)$$

since
$$\sum_{j=0}^{\infty} e^{2(j-c_X)}\left(1-\frac{\lambda}{\mu}\right)^\beta \frac{(\beta)_j}{j!}\left(\frac{\lambda}{\mu}\right)^j < \infty.$$

We also note that $G(u)$ has the expansion
$$G(u) = \sum_{k=0}^{\infty} a_k M_k\left(u; \beta, \frac{\lambda}{\mu}\right)\frac{1}{d_k},$$

$$a_k = \sum_{n=0}^{\infty} M_k\left(n; \beta, \frac{\lambda}{\mu}\right)\left(1-\frac{\lambda}{\mu}\right)^\beta \frac{(\beta)_n}{n!}\left(\frac{\lambda}{\mu}\right)^n e^{n-c_X},$$

$$\sum_{k=0}^{\infty}\frac{a_k^2}{d_k^2} = \sum_{j=0}^{\infty} e^{2(j-c_X)}\left(1-\frac{\lambda}{\mu}\right)^\beta \frac{(\beta)_j}{j!}\left(\frac{\lambda}{\mu}\right)^j < \infty.$$

Taking into account the orthogonality
$$\sum_{j=0}^{\infty} M_n\left(j; \beta, \frac{\lambda}{\mu}\right) M_m\left(j; \beta, \frac{\lambda}{\mu}\right)\left(1-\frac{\lambda}{\mu}\right)^\beta \frac{(\beta)_j}{j!}\left(\frac{\lambda}{\mu}\right)^j = \delta_n^m d_n^2, \qquad n, m \geq 1,$$

we obtain
$$\Lambda(t) = \sum_{k=0}^{\infty} a_k M_k\left(X(t); \beta, \frac{\lambda}{\mu}\right)\frac{1}{d_k},$$

$$R_\Lambda(t,s) = \mathrm{Cov}(\Lambda(t), \Lambda(s)) = \sum_{k=1}^{\infty}\frac{a_k^2}{d_k^2} e^{-k\mu|t|}$$

and
$$e^{-\mu|\tau|}\frac{a_1^2}{d_1^2} \leq R_\Lambda(\tau) \leq e^{-\mu|\tau|} \sum_{k\geq 1}\frac{M_k^2}{d_k^2}.$$

This completes the proof. □

**Remark 4.** For $\lambda > \mu$, the birth–death polynomials are again given in terms of the Meixner polynomials:
$$Q_n(j) = \left(\frac{\mu}{\lambda}\right)^n M_n\left(\frac{j}{\lambda-\mu} - \beta; \beta, \frac{\mu}{\lambda}\right), \qquad n = 0, 1, 2, \ldots,$$

and the stationary Markov process has the stationary Pascal-type distribution of the form

$$p_j = P(X(0) = (\lambda - \mu)(j + \beta)) = \left(1 - \frac{\mu}{\lambda}\right)^\beta \frac{(\beta)_j}{j!} \left(\frac{\mu}{\lambda}\right)^j, \qquad j = 0, 1, 2, \ldots,$$

so that the conditions of Propostion 1 are satisfied for $0 < q < \log \frac{\lambda}{\mu}$. Thus, a scenario of the log-Pascal type can be produced for the geometric process $X(t), t \geq 0$ by using again the orthogonality of the Meixner polynomials. We omit the details.

For $\lambda = \mu$, the birth–death polynomials are given in terms of Laguerre polynomials $L_n^{(\beta-1)}(\frac{x}{\lambda}), n = 0, 1, 2, \ldots$, as

$$Q_n(x) = \frac{n!}{(\beta)_n} L_n^{(\beta-1)}\left(\frac{x}{\lambda}\right), \qquad n = 0, 1, 2, \ldots.$$

These polynomials are orthogonal with respect to the gamma density

$$f(x) = \frac{1}{\lambda^\beta \Gamma(\lambda)} x^{\beta-1} e^{-x/\lambda}, \qquad x \geq 0.$$

The corresponding log-gamma scenario can be produced, in principle, similarly to that given in Anh, Leonenko and Shieh [[3], [4]], but a proper interpretation is lacking in this case.

### 4.3. Log-binomial scenario

**B‴**. Consider a mother process of the form

$$\Lambda(t) = \exp\{X(t) - c_X\}, \qquad c_X = N \log(p(e - 1) + 1), \tag{4.4}$$

where $X(t)$, $t \geq 0$ is a stationary birth–death process with finite state space $S = \{0, 1, \ldots, N\}$, marginal binomial distribution $\text{Bin}(N, p)$:

$$p_j = P(X(0) = j) = \binom{N}{j} (^N_j) p^j (1-p)^{N-j}, \qquad j = 0, 1, 2, \ldots, N, \ 0 < p < 1,$$

and rates

$$\lambda_n = (N - n)p, \qquad \mu_n = n(1 - p), \qquad 0 \leq n \leq N.$$

Note that

$$EX(t) = Np, \qquad \text{Var}\, X(t) = Np(1 - p),$$
$$R_X(t) = \text{Var}\, X(t) r_X(t), \qquad r_X(t) = e^{-t}, t \geq 0.$$

Under condition **B'''**, we obtain the following moment generating function:

$$M(\zeta) = \mathrm{E}\exp\{\zeta X(t)\} = (pe^\zeta + 1 - p)^N, \qquad \zeta \in \mathbb{R}.$$

We can formulate the following:

**Theorem 4.** *Suppose that condition* **B'''** *holds. Then, for any*

$$b > \frac{(p(e^2 - 1) + 1)^N}{(p(e - 1) + 1)^{2N}},$$

*the stochastic processes*

$$A_n(t) = \int_0^t \prod_{j=0}^n \Lambda^{(j)}(sb^j)\,\mathrm{d}s, \qquad t \in [0,1]$$

*converge in $L_2$ to the stochastic process $A(t)$, $t \in [0,1]$ as $n \to \infty$ such that, if $A(1) \in L_q$ for $q \in (0, \infty)$,*

$$\mathrm{E}A^q(t) \sim t^{K(q)},$$

*where*

$$K(q) = q\left(1 + \frac{c_X}{\log b}\right) - \frac{N}{\log b}\log(p(e^q - 1) + 1),$$

*and the Rényi function is given by*

$$T(q) = q\left(1 + \frac{c_X}{\log b}\right) - \frac{N}{\log b}\log(p(e^q - 1) + 1) - 1.$$

*Moreover,*

$$\mathrm{Var}\, A(t) \geq \int_0^t \int_0^t \exp(\mathrm{Var}\, X(t)e^{-\mu|u-v|} - 1)\,\mathrm{d}u\,\mathrm{d}v$$

$$\geq 2tNp(1-p)((1-e^{-t}) + 1 - e^{-t}(t+1)).$$

**Proof.** The condition for the inequalities (2.7) of Theorem 1 to hold follows from Proposition 1 since this is a process with a finite number of states, and thus $S_q$ is finite for all $q$. Regarding the condition (2.5) of Theorem 1, we consider the mother process

$$\Lambda(t) = G(X(t)), \qquad G(u) = \exp\{u - c_X\}$$

as a nonlinear transformation of the Markov process $X(t), t \in [0,1]$, which has marginal distribution $p_j \sim \mathrm{Bin}(N, p)$ and the following bivariate distribution expansion:

$$P(X(t) = j, X(s) = k) = p_j p_k \sum_{n=0}^\infty e^{-n|t-s|} K_n(j; N, p) K_n(k; N, p) \frac{1}{d_n},$$

$$d_n^2 = \frac{(-1)^n n!((1-p)/p)^n}{(-N)_n},$$

where $K_n(j; N, p), n = 0, 1, \ldots$, are Krawtchouk polynomials (see, e.g., Chihara [7] and Schoutens [29]). The recursive equation for the birth–death polynomial (3.4) becomes in this case the recursive equation for the Krawtchouk polynomials:

$$-j K_n(j; N, p) = b_n K_{n+1}(j; N, p) + \gamma_n K_n(j; N, p) + c_n K_{n-1}(j; N, p), \qquad n \geq 0,$$

$$K_{-1}(j; N, p) = 0, \qquad K_0(j; N, p) = 1,$$

and

$$b_n = p(N - n), \qquad \gamma_n = -(p(N - n) + n(1 - p)), \qquad c_n = n(1 - p).$$

Note that

$$G(u) \in l_2(\mathbb{N}, \mathrm{Bin}(N, p)),$$

since

$$\sum_{j=0}^{N} e^{2(j - c_X)} \binom{N}{j} p^j (1-p)^{N-j} < \infty.$$

We also note that $G(u)$ has the expansion

$$G(u) = \sum_{k=0}^{\infty} a_k K_n(j; N, p) \frac{1}{d_k},$$

$$a_k = \sum_{n=0}^{\infty} K_n(j; N, p) \binom{N}{j} p^j (1-p)^{N-j} e^{n - c_X},$$

$$\sum_{k=0}^{\infty} \frac{a_k^2}{d_k^2} = \sum_{j=0}^{N} e^{2(j - c_X)} \binom{N}{j} p^j (1-p)^{N-j} < \infty.$$

Taking into account the orthogonality

$$\sum_{j=0}^{N} K_n(j; N, p) K_m(j; N, p) \binom{N}{j} p^j (1-p)^{N-j} = \delta_n^m d_n^2, \qquad n, m \geq 1,$$

we obtain

$$\Lambda(t) = \sum_{k=0}^{\infty} \frac{a_k}{d_k} K_k(X(t); N, p),$$

$$R_\Lambda(t, s) = \mathrm{Cov}(\Lambda(t), \Lambda(s)) = \sum_{k=1}^{\infty} \frac{a_k^2}{d_k^2} e^{-k|t|}$$

and

$$\mathrm{e}^{-|\tau|}\frac{a_1^2}{d_1^2} \leq R_\Lambda(\tau) \leq \mathrm{e}^{-|\tau|}\sum_{k\geq 1}^{\infty}\frac{K_k^2}{d_k^2}.$$

This completes the proof. □

### 4.4. Log-hypergeometric scenario

This subsection considers a birth–death process with quadratic rates. The hypergeometric function is meant to have the power series

$$_2F_1(a,b;c;z) = \sum_{k=0}^{\infty}\frac{(a)_k(b)_k}{(c)_k}\frac{z^k}{k!} = 1 + \frac{ab}{c}z + \frac{a(a+1)b(b+1)}{c(c+1)}\frac{z^2}{2} + \cdots,$$

where $z$ is a complex variable and $a, b, c$ are parameters that can take arbitrary real or complex values provided that $c \neq 0, -1, -2, \ldots$ (see Abramowitz and Stegun [1] for details).

**B''''.** Consider a mother process of the form

$$\Lambda(t) = \exp\{X(t) - c_X\}, \qquad c_X = \log\frac{\binom{h}{N}}{\binom{g+h}{N}}{}_2F_1(-N, -g; h - N; \mathrm{e}), \qquad (4.5)$$

where $X(t)$, $t \geq 0$ is a stationary birth–death process with finite state space $S = \{0, 1, \ldots, N\}$, marginal hypergeometric distribution $\mathrm{Hyp}(N, g, h)$:

$$p_j = P(X(t) = j) = \frac{\binom{g}{j}\binom{h}{N-j}}{\binom{a+b}{N}}, \qquad j = 0, \ldots, N, \ g, h \geq N,$$

and rates

$$\lambda_n = (N-n)(g-n), \qquad \mu_n = n(h - (N-n)), \qquad n = 0, 1, \ldots, N.$$

Note that

$$\mathrm{E}X(t) = \frac{Ng}{g+h}, \qquad \mathrm{Var}\,X(t) = \frac{N(g/(g+h))(1 - g/(g+h))(g+h-N)}{a+b-1},$$

$$R_X(t) = \mathrm{Var}\,X(t)r_X(t), \qquad r_X(t) = \mathrm{e}^{-t}, \ t \geq 0.$$

Under condition **B''''**, we obtain the following moment generating function:

$$M(\zeta) = \mathrm{E}\exp\{\zeta X(t)\} = \frac{\binom{g}{N}}{\binom{g+h}{N}}{}_2F_1(-N, -g; h-N; \mathrm{e}^\zeta), \qquad \zeta \in \mathbb{R}.$$

We can formulate the following:

**Theorem 5.** *Suppose that condition* **B″″** *holds. Then, for any*

$$b > \frac{{}_2F_1(-N, -g; h-N; \mathrm{e}^2)}{\binom{h}{N}/\binom{g+h}{N}[{}_2F_1(-N, -g; h-N; \mathrm{e}^2)]^2},$$

*the stochastic processes*

$$A_n(t) = \int_0^t \prod_{j=0}^n \Lambda^{(j)}(sb^j)\,\mathrm{d}s, \qquad t \in [0,1],$$

*converge in $L_2$ to the stochastic process $A(t)$, $t \in [0,1]$ as $n \to \infty$ such that, if $A(1) \in L_q$ for $q \in (0, \infty)$,*

$$\mathrm{E}A^q(t) \sim t^{K(q)},$$

*where*

$$K(q) = q\left(1 + \frac{c_X}{\log b}\right) - \frac{N}{\log b}\log({}_2F_1(-N, -g; h-N; \mathrm{e}^q)) - \log\frac{\binom{h}{N}}{\binom{g+h}{N}},$$

*and the Rényi function is given by*

$$T(q) = q\left(1 + \frac{c_X}{\log b}\right) - \frac{N}{\log b}\log({}_2F_1(-N, -g; h-N; \mathrm{e}^q)) - \log\frac{\binom{h}{N}}{\binom{g+h}{N}} - 1.$$

*Moreover,*

$$\mathrm{Var}\,A(t) \geq \int_0^t \int_0^t \exp(\mathrm{Var}\,X(t)\mathrm{e}^{-|u-v|} - 1)\,\mathrm{d}u\,\mathrm{d}v$$

$$\geq 2t\frac{N(g/g+h)(1-g/(g+h))(g+h-N)}{g+h-1}$$

$$\times \frac{(1 - \mathrm{e}^{-((g+h+1)/2)t} + 1 - \mathrm{e}^{-((g+h+1)/2)t}(((g+h+1)/2)t+1)}{((g+h+1)/2)^2}.$$

**Proof.** The condition for the inequalities (2.7) of Theorem 1 to hold follows from Proposition 1 since this is a process with a finite number of states, and thus $S_q$ is finite for all $q$. In order to show that the condition (2.5) of Theorem 1 holds, we consider the mother process

$$\Lambda(t) = G(X(t)), \qquad G(u) = \exp\{u - c_X\}$$

as a nonlinear transformation of the Markov process $X(t), t \in [0,1]$, which has marginal distribution $p_j \sim \mathrm{Hyp}(g, h, N)$ and the following bivariate distribution expansion:

$$P(X(t) = j, X(s) = k) = p_j p_k \sum_{n=0}^\infty \mathrm{e}^{-\lambda(n)|t-s|} R_i(\lambda(n); g, h, N,) R_j(\lambda(n); g, h, N)\rho_n,$$

$$\lambda(n) = n(n-g-h-1),$$

$$\rho_n = \frac{\binom{N-h-1}{N}(-N)_n N!(-g)_n(2n-g-h-1)}{(-1)^n n!(-h)_n (n-g-h-1)_{N+1}},$$

where $R_i$ are the dual Hahn polynomials defined by

$$R_i(\lambda(n)) = R_i(\lambda(n); g, h, N) = Q_n(i, g, h, N) = Q_n(i),$$

where $Q_n(i, g, h, N), n = 0, 1, 2, \ldots$, are Hahn polynomials defined by the following recurrent relation:

$$-xQ_n(x) = hQ_{n+1}(x) - (b_n + c_n)Q_n(x) + c_n Q_{n-1}(x), \qquad Q_{-1}(x) = 0, \qquad Q_0(x) = 1,$$

where

$$b_n = \frac{(n+g+h+1)(n+g+1)(N-n)}{(2n+g+h+2)(2n+g+h+1)}, \qquad c_n = \frac{n(n+h)(n+g+h+N+1)}{(2n+g+h)(2n+g+h+1)}$$

(see, e.g., Chihara [7] and Schoutens [29]). Note that

$$G(u) \in l_2(\mathbb{N}, \text{Hyp}(g, h, N)).$$

The rest of the proof is similar to the proofs of the previous theorems using the orthogonality of Hahn polynomials and expansion of the function $G(X(t)) = e^{X(t) - c_X}$ into a series of Hahn polynomials with coefficients $a_n, n = 0, 1, 2, \ldots$, but in this case

$$\max_n \lambda(n) = \max_n \{n(n-g-h-1)\}$$

is achieved for

$$n_{\max} = \left(\frac{g+h+1}{2}\right),$$

and

$$e^{-|\tau|}\frac{a_1^2}{d_1^2} \leq R_\Lambda(\tau) \leq e^{-|\tau|n_{\max}} \sum_{k=1}^{\infty} \frac{a_k^2}{d_k^2},$$

for $d_k = 1/\sqrt{\rho_k}$. This completes the proof. $\square$

## Acknowledgements

Partially supported by the Australian Research Council Grant DP0559807, the National Science Foundation Grant DMS-04-17676, the EPSRC Grant RCMT119 and the Taiwan NSC Grant 962115M002005MY3. Leonenko's research was partially supported by the Welsh Institute of Mathematics and Computational Sciences. The authors wish to thank the referee for pointing out some obscurities in the previous version and for many constructive suggestions.